\magnification=1200
\input amssym.def
\input amssym
\input epsf
\def\CC {{\Bbb C}}
\def\DD {{\Bbb D}}

\def\PP {{\Bbb P}}
\def\QQ {{\Bbb Q}}
\def\RR {{\Bbb R}}
\def\ZZ {{\Bbb Z}}
\def\contin {\subseteq}

\def\Bar{\overline}

\def\cross {\times}
\def\mapsim {{\buildrel \sim \over \To}}
\def\part#1#2 {{\partial {#1}/\partial {#2}}}

\def\pr {\mathop{\rm pr}\nolimits}
\def\NS {\mathop{\rm NS}\nolimits}
\def\Sing {\mathop{\rm Sing}\nolimits}
\def\Spec {\mathop{\rm Spec}\nolimits}

\def\Hom {\mathop{\rm Hom}\nolimits}

\def\Pic {\mathop{\rm Pic}\nolimits}
\def\codim {\mathop{\rm codim}\nolimits}
\def\div {\mathop{\rm div}\nolimits}

\def\orb {\mathop{\rm orb}\nolimits}

\def\cE {{\cal E}}
\def\cL {{\cal L}}
\def\cN {{\cal N}}
\def\cO {{\cal O}}

\def\bda {{\bf a}}
\def\bdb {{\bf b}}

\def\bde {{\bf e}}
\def\bdf {{\bf f}}

\def\To{\longrightarrow}

\def\Sum{\sum\limits}

\def\tens{\otimes}

\def\pf{\noindent{\sl Proof: }}

\def \hf {{{1}\over{2}}}
\def\rat{\mathrel{{\hbox{\kern2pt\vrule height2.45pt depth-2.15pt
 width2pt}\kern1pt {\vrule height2.45pt depth-2.15pt width2pt}
  \kern1pt{\vrule height2.45pt depth-2.15pt width1.7pt\kern-1.7pt}
   {\raise1.4pt\hbox{$\scriptscriptstyle\succ$}}\kern1pt}}}
\def\qed{\vrule width5pt height5pt depth0pt\par\smallskip}
\def\surj{\to\kern-8pt\to}
\outer\def\startsection#1\par{\vskip0pt
 plus.3\vsize\penalty-100\vskip0pt
  plus-.3\vsize\bigskip\vskip\parskip\message{#1}
   \leftline{\bf#1}\nobreak\smallskip\noindent}
\outer\def\camstartsection#1\par{\vskip0pt
 plus.3\vsize\penalty-150\vskip0pt
  plus-.3\vsize\bigskip\vskip\parskip
   \centerline{\it#1}\nobreak\smallskip\noindent}
\def\chain{\dot{\hbox{\kern0.3em}}}
\def\cochain{\d{\hbox{\kern0.3em}}}

\def\imic{\cong}
\outer\def\thm #1 #2\par{\medbreak
  \noindent{\bf Theorem~#1.\enspace}{\sl#2}\par
   \ifdim\lastskip<\medskipamount \removelastskip\penalty55\medskip\fi}
\outer\def\prop #1 #2\par{\medbreak
  \noindent{\bf Proposition~#1.\enspace}{\sl#2}\par
   \ifdim\lastskip<\medskipamount \removelastskip\penalty55\medskip\fi}
\outer\def\lemma #1 #2\par{\medbreak
  \noindent{\bf Lemma~#1.\enspace}{\sl#2}\par
   \ifdim\lastskip<\medskipamount \removelastskip\penalty55\medskip\fi}
\outer\def\definition #1 #2\par{\medbreak
  \noindent{\bf Definition~#1.\enspace}{\sl#2}\par
   \ifdim\lastskip<\medskipamount \removelastskip\penalty55\medskip\fi}
\outer\def\corollary #1 #2\par{\medbreak
  \noindent{\bf Corollary~#1.\enspace}{\sl#2}\par
   \ifdim\lastskip<\medskipamount \removelastskip\penalty55\medskip\fi}
\def\deep #1 {_{\lower5pt\hbox{$#1$}}}
\def \restr #1 {{\big\vert_{#1}}}
\font\headfont=cmb10 scaled\magstep2

\baselineskip=15pt
\raggedbottom

\centerline{\headfont Abelian surfaces in toric 4-folds}
\bigskip
\centerline{G.K. Sankaran}
\medskip
\centerline{Department of Mathematical Sciences, University of Bath}
\centerline{Bath BA2 7AY, England}
\centerline{\tt gks@maths.bath.ac.uk}

\noindent Mathematics Subject Classification: 14E25, 14K05, 14M25, 14F05

\bigskip
There are embeddings of complex abelian surfaces in $\PP^4$ but it was
shown by Van de Ven in [17] that no abelian $d$-fold can be embedded
in $\PP^{2d}$ if $d\ge 3$. Hulek [9], Lange [13], Birkenhake [4] and
Bauer and Szemberg [3] have all investigated the possibility of
replacing $\PP^{2d}$ with a product of projective spaces. Furthermore,
Lange [14] studies the rank~$2$ bundle on $\PP^1\cross\PP^3$ that
arises from the abelian surfaces in $\PP^1\cross\PP^3$ by Serre's
construction. The analogous bundle associated with the abelian
surfaces in $\PP^4$ is, of course, the Horrocks-Mumford bundle.

In this paper we shall work over the complex numbers and consider
embeddings of abelian surfaces in slightly more general ambient spaces
of dimension~$4$, namely smooth toric varieties. The most tractable
and probably the most interesting cases are when the ambient
variety has small Picard number. We shall therefore consider the
following question. Suppose $X$ is a smooth complete toric variety of
dimension~$4$ and $\rho(X)\le 2$. Does there exist an abelian surface
$A\contin X$ and, if so, can we describe the embedding?

In the first section we shall give some numerical conditions that such
an embedding must satisfy and show that for many $X$ of this type
there can be no totally nondegenerate (see Definition~1.1, below)
abelian surfaces in~$X$. In Section~2, which is based on unpublished
joint work with Professor T.~Oda, we show how to describe morphisms
into smooth toric varieties in a particularly simple way. The results
of this section overlap with work of Cox [6], Guest [8] and Kajiwara
[11] but it is useful to us to have them in the form given here. We
apply this description in Section~3, in which we exhibit a new
$2$-dimensional family of abelian surfaces embedded in a particular
toric $4$-fold $X$ of Picard number~$2$.

The normal bundles of the surfaces described in Section~3 give rise to
rank~$2$ bundles on~$X$ which should be interesting to study. However,
we do not attempt this here, but in Section~$4$ we make a few comments
on this and other related matters.

Much of this work was carried out in 1997 while the author was
visiting the Research Institute for Mathematical Sciences in Kyoto,
supported by the JSPS International Project Research ``Infinite
Analysis and Geometry''.

\startsection 1. Numerical conditions.

If $X$ is a smooth toric $4$-fold and $\rho(X)=1$ then $X\imic\PP^4$
and the only possibility is that $A$ is a Horrocks-Mumford
surface. So we consider the case $\rho(X)=2$.

\definition 1.1 An abelian surface $A\contin X$, $X$ a toric $4$-fold, is
{\it totally nondegenerate}\/ if $A\cap X_\sigma$ is of dimension~$1$
for every torus-invariant divisor $X_\sigma\contin X$.

We shall be interested only in totally nondegenerate embeddings. An
example of an embedding that fails to be totally nondegenerate may be
obtained by taking a Horrocks-Mumford surface $A\contin\PP^4$ and
taking $X$ to be the blow-up of $\PP^4$ in a point outside~$A$.

Smooth toric $4$-folds of Picard number $2$ are well understood. In
fact, smooth complete toric varieties of Picard number $2$ in any
dimension have been classified by Kleinschmidt [12]. Such a toric
variety is a projectivisation of a decomposible vector bundle over a
projective space of smaller dimension. So in our case $X$ is a
$\PP^3$-bundle over $\PP^1$, a $\PP^2$-bundle over $\PP^2$, or a
$\PP^1$-bundle over $\PP^3$.

\thm 1.2 If $X$ is the projectivisation of a decomposible
$\PP^3$-bundle over $\PP^1$ then $X$ contains no totally nondegenerate
abelian surfaces unless $X=\PP^1\cross\PP^3$.

The proof will be given as part of the analysis below. First, we want
a toric description of $X$ (this is a convenient way to do the
calculations).  We can write 
$$
X=\PP\big(\cO\oplus\cO(\kappa_1)
\oplus\cO(\kappa_2)\oplus\cO(\kappa_3)\big)
$$
and without loss of generality we may suppose that
$\kappa_1\ge\kappa_2 \ge\kappa_3\ge 0$.  We put
$\kappa=\kappa_1+\kappa_2+\kappa_3$. Then $X=X_\Sigma$, where $\Sigma$
is the fan in $\RR^4$ whose $1$-skeleton consists of
$\sigma_1=(1,0,0,0)$, $\sigma_2=(0,1,0,0)$ $\sigma_3=(0,0,1,0)$,
$\sigma_4=(-1,-1,-1,0)$ (these four form a primitive collection in the
sense of Batyrev [1]), $\tau_1=(0,0,0,1)$ and
$\tau_2=(\kappa_1,\kappa_2,\kappa_3,-1)$, and whose top-dimensional
cones are spanned by a $\tau$ and three of the $\sigma$s.

Let $D_i=\Bar{\orb\sigma_i}$ and $E_j=\Bar{\orb\tau_j}$. Then $\Pic X
=\ZZ\bda\oplus\ZZ\bdb$, where $\bda=[E_1]$ and $\bdb=[D_4]$. Note that
$\bda=[E_1] =[E_2]$ is the class of a fibre of the projection
$p:X\to\PP^1$ and that the restriction of $\bdb$ to a fibre is
$\cO_{\PP^3}(1)$. Also $[D_i]-[D_4]=-\kappa_i[E_2]$ so
$[D_i]=\bdb-\kappa_i\bda$ for $i=1,\ 2,\ 3$. The intersection numbers
are $\bda^4=\bda^3\bdb=\bda^2\bdb^2=0$ (in fact $\bda^2=0$ in
$H^4(X;\ZZ)$ since $\bda$ is the class of a fibre), $\bda\bdb^3=1$ and
$\bdb^4=\kappa$.

Now suppose $A\contin X$ is an abelian surface. The class of $A$ in
$H^4(X;\ZZ)$ (or in $A^2(X)$, which is the same thing in this case) is
$\lambda\bda\bdb+\mu\bdb^2$ for some $\lambda,\mu\in\ZZ$. By
Proposition~3 of [17], which is a version of the self-intersection
formula in [7]
$$
c_2(\cN_{A/X})\cdot[A]=[A]^2.
$$
We have $c(\cN_{A/X})=c(X)/c(A)$ and $c(A)=1$ so
$c_2(\cN_{A/X})=c_2(X)$. The total Chern class of a smooth complete
toric variety is well known ([16], Theorem~3.12). Here
$$
c(X)=\prod_{i=1}^4\Big(1+[D_i]\Big)\prod_{j=1}^2\Big(1+[E_j]\Big)
$$
and in degree~2
$$
\eqalign{
c_2(X)&=\Sum_{i<j}[D_i][D_j]+\Sum_{i,j}[D_i][E_j]\cr
      &=\Sum_{i<j}(\bdb-\kappa_i\bda)(\bdb-\kappa_j\bda)+8\bda\bdb\cr
      &=(8-3\kappa)\bda\bdb+6\bdb^2.\cr}
$$
So, by the self-intersection formula,
$$
\big((8-3\kappa)\bda\bdb+6\bdb^2\big)\big(\lambda\bda\bdb+\mu\bdb^2\big)
=\big(\lambda\bda\bdb+\mu\bdb^2\big)^2
$$
which simplifies to
$$
\lambda(2\mu-6)=\kappa(3\mu-\mu^2)+8\mu.\eqno{\hbox{($\dagger$)}}
$$
Next, note that $\mu=[A]\bda\bdb$ which is the degree of the space
curve obtained by intersecting $A$ with a general fibre of
$p:X\to\PP^1$. As this curve is contained in an abelian variety it
cannot be rational, so $\mu\ge 3$. Put $\nu=\mu-3$. From ($\dagger$)
we know that $\nu\ne 0$ and
$$
\eqalign{
2\lambda\nu&=\kappa\big(3\nu+9-(\nu+3)^2\big)+8\nu+24\cr
           &=-\kappa(\nu^2+3\nu)+8\nu+24\cr}
$$
so $2\nu|-\kappa(\nu^2+3\nu)+24$.

Put $\nu=2^r\nu'$ with $\nu'$ odd. Then
$$
2^{r+1}\nu'|-\kappa(2^{2r}\nu'^2+3\cdot 2^r\nu')+24
$$
so $\nu'|24$ so $\nu'=1$ or $\nu'=3$. Moreover, if $\nu'=1$ we have
$$
2^{r+1}|-\kappa(2^{2r}+3\cdot 2^r)+24
$$
so $r\le 2$ or $2^{r-2}|-3\cdot 2^{r-3}\kappa+3$ and so $r=3$. Thus if
$\nu'=1$ then $\nu=1$, $2$, $4$ or $8$ and $\mu=4$, $5$, $7$ or
$11$. Similarly, if $\nu'=3$, then
$$
2^{r+1}|-\kappa(9\cdot 2^{2r}+9\cdot 2^r)+24
$$
so $r\le 2$ or $2^{r-2}|-9\cdot 2^{r-3}\kappa+3$ and again $r=3$. So
if $\nu'=3$ then $\nu=3$, $6$, $12$ or $24$ and $\mu=6$, $9$, $15$ or $27$.

Consider the curves $B_i=A\cap D_i$ on~$A$. We have
$$
p_g(B_i)=\hf[A](\bdb-\kappa_i\bda)^2+1
        ={{\lambda+\kappa\mu}\over{2}}-\kappa_i\mu+1
$$
and this must of course be a positive integer. Adding together the
inequalities
$$
0\le
{{\lambda}\over{2}}+\big({{\kappa}\over{2}}-\kappa_i\big)\mu
\eqno{\hbox{($*_i$)}}
$$
for $i=1$, $2$, $3$ and using ($\dagger$),
$$
0\le 3\lambda+\kappa\mu = {{12\mu}\over{\mu-3}}-{{\kappa\mu}\over{2}}.
$$
If equality holds here then
$0={{\lambda}\over{2}}+\big({{\kappa}\over{2}}-\kappa_i\big)\mu$ for
all~$i$, so $\kappa_1=\kappa_2=\kappa_3$ and $3|\kappa$. So
$$
\kappa\le{{24}\over{\mu-3}}={{24}\over{\nu}}
$$
with equality only if $3|\kappa$.

For $\mu=27$ this implies $\kappa=0$ and then $X=\PP^3\cross\PP^1$
which is treated in~[9]. In any case it does not occur as it gives
$p_g(B_4)={{9}\over{4}}$. If $\mu=15$ then $\kappa=0$ or $1$ and in
both cases $p_g(B_4)$ fails to be an integer. Similarly if $\mu=9$ we
have $\kappa\le 3$ but $p_g(B_4)=3+{{9}\over{4}}\kappa$ so $4|\kappa$,
so $\kappa=0$ and $X=\PP^3\cross\PP^1$ (and according to [9] this
case does not occur either). We shall see shortly that $\mu=11$, $7$,
$5$ or $4$ is impossible for a different reason, so we are left with
$\mu=6$.

If $\mu=6$ then $\kappa<8$ and $p_g(B_4)=4+{{3}\over{2}}\kappa$ so
$\kappa=0$, $2$, $4$ or $6$. The case $\kappa=0$ is covered by [9]
and [13] (and this case really does occur). If $\kappa=4$ then
$\kappa\ge 2$ and then $\lambda=8-3\kappa=-4$ so $*_1$
fails. Similarly if $\kappa=6$ then $\lambda=-10$ so $*_1$ fails
unless $\kappa_1\le 2$, which implies
$\kappa_1=\kappa_2=\kappa_3=2$. If $\kappa=\kappa_1=2$ then $*_1$
fails so the remaining case is $\kappa=2$, $\kappa_1=\kappa_2=1$,
$\kappa_3=0$.

However, neither of these cases is possible, because in either case
$h^0(\cO_A(B_1))=1$ since $B_1^2=2$, so $|B_1|$ is a point and
therefore $B_1=B_2$. But then $A$ is contained in the closure of a
smaller torus, namely $\{(t_1,t_1,t_3,t_4)\}\imic(\CC^*)^3$, and no
abelian surface can be embedded in a smooth toric 3-fold.

It remains to eliminate the possibilities $\mu=4,\ 5,\ 7,\ 11$. By a
standard theorem ([15], Section~3.3, or [4]) there is a commutative
diagram
$$
\matrix{0&\To&C&\To&A&\To&C'&\To&0\cr
         &&&&&\searrow^p&\downarrow&&\cr
         &&&&&&\PP^1&&\cr}
$$
where the top row is an exact sequence of abelian varieties, so $C$
and $C'$ are elliptic curves. A general fibre $F$ of $p:A\to\PP^1$
is therefore a disjoint union of $d$~translates of~$C$, where $d$ is
the degree of $C'\to\PP^1$ (and therefore $d\ge 2$). Now
$\mu=[A]\bda\bdb=(F.B_4)_A=d(C.B_4)_A$, and since $C.B_4$ is also
the degree of $C\contin\PP^3$ we have $C.B_4\ge 3$. This shows
that $\mu$ cannot be equal to $4$, $5$, $7$ or $11$ and completes the
proof of the theorem.~\qed

Next we consider the case where $X$ is a $\PP^1$-bundle over $\PP^3$,
which is easy.

\thm 1.3 If $X$ is a $\PP^1$-bundle over $\PP^3$ then $X$ contains no
totally nondegenerate abelian surfaces unless $X=\PP^1\cross\PP^3$.

\pf Suppose $X=\PP\big(\cO\oplus\cO(\kappa)\big)$, $\kappa>0$. Let
$p:X\to\PP^3$ be the projection and put
$\bda=[p^*\cO_{\PP^3}(1)]$. Take coordinates $(x,y)$ on
$\cO\oplus\cO(\kappa)$ and put $\bdb=[(x=0)]$. Then
$\bdb-\kappa\bda=[(y=0)]$ is also the class of a section, and since
$(y=0)$ is disjoint from $(x=0)$ we have $\bdb(\bdb-\kappa\bda)=0$ in
$H^4(X;\ZZ)$. So $H^4(X;\ZZ)$ is generated by $\bda^2$ and
$\bda\bdb$. Suppose $A$ is an abelian surface in $X$ and that
$[A]=\lambda\bda^2+\mu\bda\bdb$. Since $\bda|_A $ and
$(\bdb-\kappa\bda)|_A $ are disjoint effective curves on $A$
neither can be ample, but any effective class with positive
self-intersection on an abelian surface is ample (see [15]). So
$\bda^2[A]=\bda(\bdb-\kappa\bda)[A]=(\bdb-\kappa\bda)^2[A]=0$, which,
combined with the intersection numbers $\bda^4=0$, $\bda^3\bdb=1$,
$\bda^2\bdb^2=\kappa$, $\bda\bdb^3=\kappa^2$ and $\bdb^4=\kappa^3$,
gives $\lambda=\mu=0$ if $\kappa\ne 0$. This is impossible.~\qed

For the remaining case, when $X$ is a $\PP^2$-bundle over $\PP^2$, the
methods above do not suffice to determine a finite list of possible
cases. However, we can give some quite strong necessary
conditions. Suppose then that
$$
X=\PP\big(\cO_{\PP^2}\oplus\cO_{\PP^2}(\kappa_1)
\oplus\cO_{\PP^2}(\kappa_2)\big)
$$
with $\kappa_1\ge \kappa_2 \ge 0$.
Put $\kappa=\kappa_1+\kappa_2$. Then $X=X_\Sigma$, where the
$1$-skeleton of the fan $\Sigma$ consists of 
$$\displaylines{
\sigma_1=(1,0,0,0),\qquad\sigma_2=(0,1,0,0),\qquad\sigma_3=(-1,-1,0,0),\cr
\tau_1=(0,0,1,0),\qquad\tau_2=(0,0,0,1),\qquad\tau_3=
(\kappa_1,\kappa_2,-1,-1),\cr
}
$$
and the top-dimensional cones are spanned by two $\sigma$s and two
$\tau$s. If we put $D_i=\Bar{\orb\sigma_i}$ and $E_j=\Bar{\orb
\tau_j}$ and $\bda=[E_1]$, $\bdb=[D_3]$, then
$[E_1]=[E_2]=[E_3]=\bda=p^*\cO_{\PP^2}(1)$, where $p:X\to\PP^2$ is the
projection, and $[D_i]=\bdb-\kappa_i\bda$. The intersection numbers
are $\bda^4=\bda^3\bdb=0$, $\bda^2\bdb^2=1$, $\bda\bdb^3=\kappa$ and
$\bdb^4=\kappa^2$.

If $A\contin X$ is an abelian surface we can take
$[A]=\lambda'\bda^2+\mu'\bda\bdb+\nu'\bdb^2\in H^4(X;\ZZ)$. The
notation is convenient because it is easier to work with $\nu=\nu'$,
$\mu=\mu'+\kappa\nu$ and $\lambda=\lambda'+\kappa\mu$ than with
$\lambda'$, $\mu'$ and $\nu'$ directly. Then $\bda^2[A]=\nu$,
$\bda\bdb[A]=\mu$ and $\bdb^2[A]=\lambda$, so $\lambda$, $\mu$ and
$\nu$ are all non-negative and $\lambda$ and $\nu$ are even. We assume
that $\kappa>0$, since otherwise $X=\PP^2\cross\PP^2$ and then by [9]
we know that $A$ is the product of two plane cubics. We also make the
nondegeneracy assumption that $\nu>0$, that is, that $p:A\to\PP^2$ is
surjective. Now the Hodge index theorem on $A$ gives
$$
\lambda\nu\le\mu^2\eqno{(1)}
$$
and the self-intersection formula gives
$$
(3-3\kappa+\kappa_1\kappa_2)\nu+(9-2\kappa)\mu+3\lambda
=2\lambda\nu-2\kappa\mu\nu+\mu^2\eqno{(2)}
$$

Since $D_i|_A\ge 0$ we also have $(\bdb-\kappa_1\bda)^2[A]\ge 0$ and
$\bda(\bdb-\kappa_1\bda)[A]\ge 0$ so
$$
\lambda-2\kappa_1\mu+\kappa_1^2\nu\ge 0\eqno{(3)}
$$
and
$$
\mu-\kappa_1\nu\ge 0\eqno{(4)}
$$

We can rewrite all of these in terms of $x=\mu/\kappa\nu$ and
$y=\lambda/\kappa^2\nu$:
$$
\eqalign{
y&\le x^2\cr
y&={{-\nu}\over{2\nu-3}}x^2+{{2\nu-2+9/\kappa}\over{2\nu-3}}x
+{{3+\kappa_1\kappa_2-3\kappa}\over{\kappa^2(2\nu-3)}} =f(x)\cr
y&\ge{{2\kappa_1}\over{\kappa}}x-{{\kappa_1^2}\over{\kappa^2}}\cr
x&\ge{{\kappa_1}\over{\kappa}}\cr}
$$
The three inequalities are satisfied for $(x,y)$ in the shaded area in
the diagram.

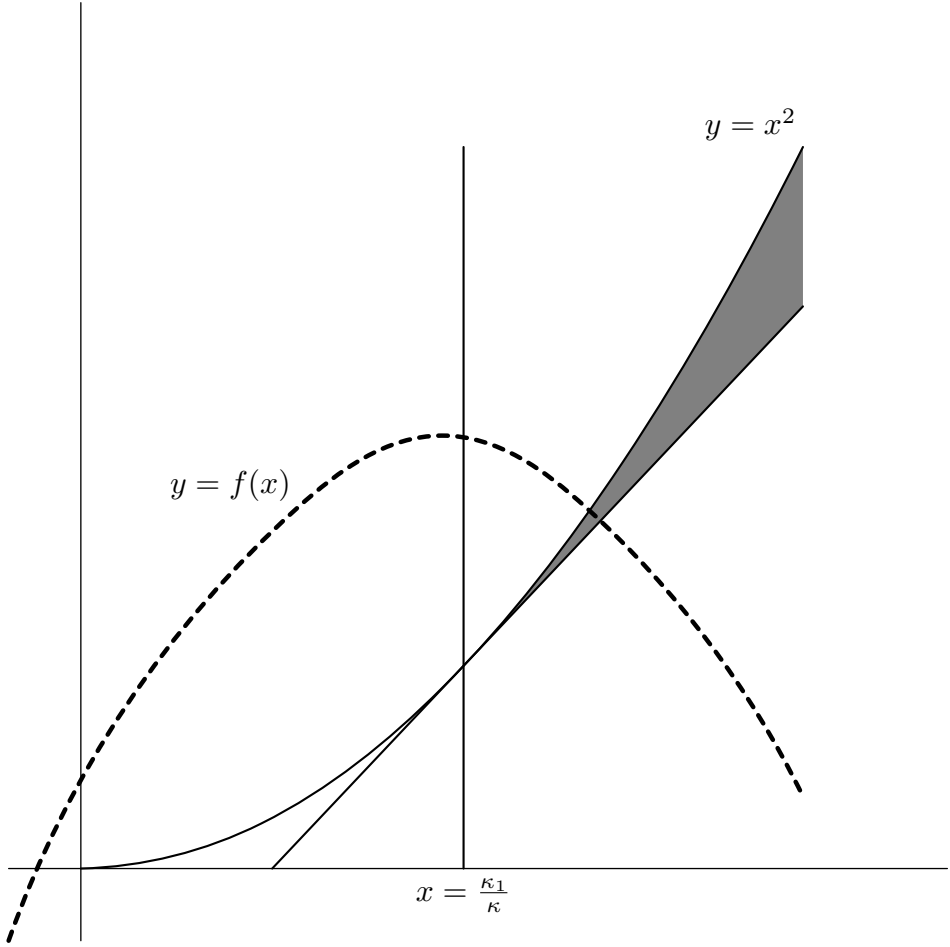
\epsffile{embedfig.ps}

From this we can deduce the following (tidy but not very sharp)
result.

\thm 1.4 Suppose
$X=\PP\big(\cO_{\PP^2}\oplus\cO_{\PP^2}(\kappa_1)
\oplus\cO_{\PP^2}(\kappa_2)\big)$ with $\kappa_1\ge \kappa_2 \ge 0$
and $p:X\to\PP^2$ is the projection. If $\kappa_1>2\kappa_2$ then $X$
contains no totally nondegenerate abelian surface $A$ for which
$p:A\to\PP^2$ is surjective.

\pf Clearly the curve with equation $y=f(x)$ will not pass through the
shaded area if we have $f(\kappa_1/\kappa)<\kappa_1^2/\kappa^2$ and
$f'(\kappa_1/\kappa)\le 2\kappa_1/\kappa$. So if these inequalities
hold no such abelian surface will exist.  Moreover, if such a surface
does exist then $\nu\ge 6$ since by Riemann-Roch
$$
h^0\big(p^*\cO_{\PP^2}(1)\big)=\hf\bda^2[A]=\hf\nu
$$
and $h^0\big(p^*\cO_{\PP^2}(1)\big)\ge h^0\big(\cO_{\PP^2}(1)\big)=3$.

I claim that in fact $f'(\kappa_1/\kappa)\le 2\kappa_1/\kappa$ unless
$\kappa=1$, when $\kappa_1=1> 2\kappa_2=0$: we deal with this
possibility below. For if $f'(\kappa_1/\kappa)>2\kappa_1/\kappa$ then
$$
{{-2\nu}\over{2\nu-3}}\,{{\kappa_1}\over{\kappa}}
+{{2\nu-2+9/\kappa}\over{2\nu-3}}> 2{{\kappa_1}\over{\kappa}}
$$
so
$$
\eqalign{
(2\nu-2)\kappa+9&> (6\nu-6)\kappa_1\cr
                &\ge (3\nu-3)\kappa\cr}
$$
since $\kappa_1\ge \kappa_2=\kappa-\kappa_1$. So $9\ge(\nu-1)\kappa\ge
10$ unless $\kappa=1$.

Thus no abelian surface as in the theorem will exist if
$f(\kappa_1/\kappa)<\kappa_1^2/\kappa^2$, that is, if
$$
-{{\nu}\over{2\nu-3}}\,{{\kappa_1^2}\over{\kappa^2}}
+{{2\nu-2+9/\kappa}\over{2\nu-3}}\,{{\kappa_1}\over{\kappa}}
+{{3+\kappa_1\kappa_2-3\kappa}\over{\kappa^2(2\nu-3)}}
<{{\kappa_1^2}\over{\kappa^2}}
$$
which simplifies to
$$
\kappa_1^2-{{2\nu-1}\over{\nu+1}}\kappa_1\kappa_2
-{{6}\over{\nu+1}}\kappa_1 +{{3}\over{\nu+1}}(\kappa_2-1) >0.
$$

We first assume $\kappa_2\ge 1$. Then this will certainly hold if
$$
\kappa_1^2-\left(2-{{3}\over{\nu+1}}\right)\kappa_1\kappa_2
-{{6}\over{\nu+1}}\kappa_1>0,
$$
that is, if
$$
\kappa_1>\left(2-{{3}\over{\nu+1}}\right)\kappa_2+{{6}\over{\nu+1}}
$$
which is true if $\kappa_1>2\kappa_2$.

It remains to deal with the possibility that $\kappa_2=0$. Then there
are no abelian surfaces as long as
$$
(\nu+1)\kappa_1^2-6\kappa_1-3>0
$$
and since $\nu\ge 6$ and is even the only possibilities are
$\kappa=\kappa_1=1$, $\nu=6$ or $\nu=8$. If $\nu=6$ then by (1) and
(2)
$$
9\lambda=-\mu^2+19\mu\ge 54
$$
so $\mu=9$ or $\mu=10$ and in both cases $\lambda=10$ and (3) fails.
If $\nu=8$ then an identical argument shows that $\lambda=10$ and
$\mu=10$ or $\mu=13$, contradicting (3).~\qed

\startsection 2. Morphisms to smooth toric varieties

This section is based on joint work with Tadao Oda. I am grateful to
Professor Oda for allowing me to use these results here.

Let $\Delta$ be a finite (but not necessarily complete) smooth fan for
a free $\ZZ$-module $N\imic\ZZ^r$ of rank $r$, and denote by $X$ and
$T:=T_N$ the corresponding toric variety and the algebraic torus. For
simplicity we work over~$\CC$. $M:=\Hom(N,\ZZ)$ is the $\ZZ$-module
dual to $N$ with the duality pairing $<\;,\;> :M\cross N\to\ZZ$. As a
general reference for the theory of toric varieties, we use [16].
 
As usual, $\Delta(1)$ denotes the set of one-dimensional cones in
$\Delta$. For each $\rho\in\Delta(1)$, we denote by $V(\rho)$ the
corresponding irreducible Weil divisor $\Bar{\orb\rho}$ on $X$.
 
\thm 2.1 Let $Y$ be a normal algebraic variety over $\CC$. Then the set of
morphisms $f: Y\to X$ such that $f(Y)\cap T\neq\emptyset$ is in
one-to-one correspondence with the set of pairs
$(\{D(\rho)\}_{\rho\in\Delta(1)},\varepsilon)$ consisting of a set
$\{D(\rho)\}_{\rho\in\Delta(1)}$ of effective Weil divisors $D(\rho)$
on $Y$ for $\rho\in\Delta(1)$ such that
$$
D(\rho_1)\cap D(\rho_2)\cap\cdots\cap D(\rho_s)=\emptyset\qquad
\hbox{whenever }\rho_1+\rho_2+\cdots+\rho_s\not\in\Delta
$$
and a group homomorphism
$$
\varepsilon:M\to H^0\Bigl(Y\setminus{\textstyle\bigcup\limits}_{\rho\in\Delta(1)}D(\rho),
\cO_Y\Bigr)^{\cross}
$$
to the multiplicative group of invertible regular functions on
$Y\setminus\bigcup\limits_{\rho\in\Delta(1)}D(\rho)$ such that
$$
\div(\varepsilon(m))=\sum_{\rho\in\Delta(1)}
< m,n(\rho)> D(\rho)\qquad\hbox{for all}\quad m\in M.
$$
 
\pf Suppose a morphism $f:Y\to X$ with $f(Y)\cap T\neq\emptyset$
is given. For each $\rho\in\Delta(1)$, the pull-back Weil divisor
$D(\rho):=f^{-1}(V(\rho))$ is well-defined, since $Y$ is assumed to be
normal and $X$ smooth and $f(Y)\not\subset V(\rho)$.
 
If $\rho_1,\ldots,\rho_s\in\Delta(1)$ satisfy
$\rho_1+\cdots+\rho_s\not\in\Delta$, then we obviously have
$V(\rho_1)\cap V(\rho_2)\cap\cdots\cap V(\rho_s)=\emptyset$, hence
$D(\rho_1)\cap D(\rho_2)\cap\cdots\cap D(\rho_s)=\emptyset$. By
assumption, $f^{-1}(T)=Y\setminus\bigcup_{\rho\in\Delta(1)}D(\rho)$ is
a nonempty open set of $Y$, and the restriction of $f$ to it induces a
ring homomorphism
$$
f^*:\CC[M]\to H^0(Y\setminus{\textstyle\bigcup\limits}_{\rho\in\Delta(1)}D(\rho),
\cO_Y),
$$
where $\CC[M]:=\bigoplus\limits_{m\in M}\CC\bde(m)$ is the semigroup ring
of $M$ over $\CC$ so that $T=\Spec(\CC[M])$. The composite
$\varepsilon:=f^*\circ\bde$ with $\bde:M\to\CC[M]$ obviously
satisfies our requirements, since
$$
\div(\bde(m))=\sum_{\rho\in\Delta(1)}< m,n(\rho)> V(\rho)
\qquad\hbox{for all}\quad m\in M.
$$
 
Conversely, suppose $(\{D(\rho)\}_{\rho\in\Delta(1)},\varepsilon)$
satisfying the requirements are given. Put
$$
\hat{\sigma}:=\{\rho\in\Delta(1)\;\mid\;\rho\not\prec\sigma\}
\qquad\hbox{for}\quad\sigma\in\Delta.
$$
Then we have
$$
\eqalign{
\Spec(\CC[M\cap\sigma^{\vee}])&=U_{\sigma}\cr
&=\bigcap_{\rho\in\hat{\sigma}}(X\setminus V(\rho))\cr
&=X\setminus\bigcup_{\rho\in\hat{\sigma}}V(\rho)\cr}.
$$
If we denote
$Y_{\sigma}:=Y\setminus\bigcup_{\rho\in\hat{\sigma}}D(\rho)$, then we
have $Y=\bigcup_{\sigma\in\Delta}Y_{\sigma}$. Indeed, the right hand
side is the complement in $Y$ of
$\bigcap_{\sigma\in\Delta}(\bigcup_{\rho\in\hat{\sigma}}D(\rho))$,
which is the union of $\bigcap_{\sigma\in\Delta}D(\rho(\sigma))$ for
all $\{\rho(\sigma)\in\Delta\;\mid\;\rho(\sigma)\in\hat{\sigma},
\forall\sigma\in\Delta\}$, hence is empty by assumption.
 
For each $\sigma\in\Delta$, $M\cap\sigma^{\vee}$ is the semigroup
consisting of $m\in M$ such that $\bde(m)$ is regular on $U_{\sigma}$.
By assumption, we thus see that $\varepsilon(M\cap\sigma^{\vee})$
consists of regular functions on $Y_{\sigma}$. Hence we get a morphism
$f_{\sigma}:Y_{\sigma}\to U_{\sigma}$. Clearly, we can glue
$\{f_{\sigma}\}_{\sigma\in\Delta}$ together to get a morphism
$f:Y\to X$.~\qed

Although we shall not need it in the rest of the paper we mention here
a simple consequence of this result and an example.

\corollary 2.2 Let $y_0\in Y$ be a point of a normal algebraic variety
$Y$ over $\CC$. Then the set of morphisms $f:Y\to X$ such that
$f(y_0)$ coincides with the identity element $1\in T$ is in one-to-one
correspondence with the set of pairs
$(\{D(\rho)\}_{\rho\in\Delta(1)},\varepsilon)$ satisfying the same
conditions as in Theorem 2.1 and such that
$\varepsilon(m)$ has value $1$ at $y_0$ for all $m\in M$.~\qed
 
Let us consider the case (first investigated by Guest in [8]) where
$Y=\PP^1$ is the projective line with $y_0=\infty$. The morphisms
$f:Y\to X$ satisfying $f(\infty)=1$ are in one-to-one correspondence
with the pairs $(\{D(\rho)\}_{\rho\in\Delta(1)},\varepsilon)$
satisfying the conditions of Corollary~2.2, so that $\varepsilon(m)$
has value $1$ at $\infty$ for all $m\in M$.
 
In terms of an inhomogeneous coordinate $z$ on $Y=\PP^1$, let us
identify the effective divisor $D(\rho)$ as usual with a monic
polynomial $P_{\rho}(z)\in\CC[z]$ for each $\rho\in\Delta(1)$. Then
for each $m\in M$ we have
$$
\varepsilon(m)=
\prod_{\rho\in\Delta(1)}P_{\rho}(z)^{< m,n(\rho)>}.
$$
Our requirements amount to the following:
$P_{\rho_1},P_{\rho_2},\ldots,P_{\rho_s}$ have no common factors
whenever $\rho_1+\rho_2+\cdots+\rho_s\not\in\Delta$, and
$$
\sum_{\rho\in\Delta(1)}< m,n(\rho)>\deg P_{\rho}=0
\qquad\hbox{for all}\quad m\in M.
$$

A strong result describing morphisms into toric varieties is proved by
Cox in [6]. The version given here is perhaps simpler to apply but is
much more limited in its scope. It has been extended by Kajiwara [11]
to certain singular toric varieties~$X$, including all projective
toric varieties.

\startsection 3. An example

In this section we shall show that one possibility not excluded by the
results of section 1 does indeed occur.

\thm 3.1 There is a $2$-dimensional family of abelian surfaces
$A\contin X=\PP\big(\cO_{\PP^2}\oplus
\cO_{\PP^2}(1)\oplus\cO_{\PP^2}(1)\big)$ such that $[A]=-6\bda^2+2\bda\bdb+6\bdb^2$.

Before starting to prove this theorem, let us compare this case with
the restrictions given in Theorem 1.4. It is really the simplest case
not excluded already. We have taken $\kappa_1=\kappa_2=1$, thus
complying with 1.4, and $\nu=\nu'=6$, which is minimal. Then the
values $\lambda=22$, $\mu=14$ are dictated by the equations
(1--4) of Section~1. In fact the inequality (3) is in this case an equality:
geometrically this means that the abelian surface $A$ which arises
turns out to be isogenous to a product of two elliptic curves.

The strategy for proving Theorem~3.1 is as follows. We first show that
there exist abelian surfaces having curves which behave numerically
like the intersections of a surface of class
$-6\bda^2+2\bda\bdb+6\bdb^2$ with toric strata in~$X$. Given such a
surface $A$, we apply Theorem~2.1 to obtain a morphism $\phi:A\to
X$. This morphism will depend on the choice of the curves. We show
also, again using Theorem~2.1, that such a choice of curves also
determines a morphism $\psi:A\to\PP^2\cross\PP^1$ and that, for a
general choice of $A$ and of the curves, $\psi$ is birational onto its
image. Furthermore, $\phi$ factors through $\psi$ as a rational map
and is therefore also birational onto its image. We describe the
singular locus of $\psi(A)$ and show that, for general $A$, we can
choose things so that $\phi(A)$ has isolated singularities. Then by an
application of the double point formula we can deduce that $\phi(A)$
is smooth.

\prop 3.2 There exists a $2$-dimensional family of abelian surfaces
$A$ containing curves $E_1$ and $C$ such that $E_1^2=6$, $C^2=0$
and $E_1.C=4$.

\pf Take $A=\CC^2/\Lambda$, where $\Lambda$ is the lattice spanned
by the columns $\bdf_i$ of the period matrix
$$
\Pi=\pmatrix{4\tau_1&3\tau_1&1&0\cr 3\tau_1&\phantom{3}\tau_3&0&3\cr}
$$
with $\pmatrix{4\tau_1&3\tau_1\cr 3\tau_1&\phantom{3}\tau_3\cr}$ in
the Siegel upper half-space of degree~$2$. The complex torus $A$ is
then an abelian surface equipped with a polarisation $H$ of type
$(1,3)$. We take $E_1$ to be some curve on $A$ giving rise to this
polarisation. Additionally, $A$ contains the elliptic curve
$\CC/\ZZ+\ZZ\tau_1$, embedded by $\gamma:z\mapsto\pmatrix{4z\cr
3z}$: we take $C$ to be this image. Then
$\gamma(1)=\pmatrix{4\cr 3}=4\bdf_3+3\bdf_4$ and
$\gamma(\tau_1)=\pmatrix{4\tau_1\cr 3\tau_1}=\bdf_1$, so
$E_1.C=\deg_H C=H(\bdf_1,4\bdf_3+3\bdf_4)=4$ as required.~\qed

According to Theorem~2.1, we must now specify the curves $E_1$, $E_2$,
$E_3$, $D_1$, $D_2$, $D_3$ on $A$ and also specify trivialisations of
certain line bundles on $A$. We choose $D_1$ to be some element of the
linear system $|2C|$. This is of dimension~$1$, because in
$H^2(A,\ZZ)\imic\bigwedge^2\Hom(\Lambda,\ZZ)$ we have
$$
[C]=\bdf^*_1\wedge 4\bdf^*_3+\bdf^*_1\wedge\bdf^*_4
$$
which is not divisible. Hence $\cO_A(C)$ is of type $(0,1)$ and
according to [15] it follows that $\dim H^0\big(\cO_A(2C)\big)=2$. A
general $D_1$ in this system is a union of two disjoint elliptic
curves, both translates of $C$. We choose $E_1$ to be an element of
the polarising class. Now choose homogeneous coordinates $(x_1:x_2)$
in $|D_1|\imic\PP^1$ and $(y_1:y_2:y_3)$ in $|E_1|\imic\PP^3$ such
that $D_1=(1:0)$ and $E_1=(1:0:0)$ and put $D_2=(0:1)\in|D_1|$,
$E_2=(0:1:0), E_3=(0:0:1)\in|E_1|$.

\thm 3.3 The complete linear systems $|D_1|$ and $|E_1|$ determine a
morphism $\psi=(\phi_{|E_1|},\phi_{|D_1|}): A\to\PP^2\cross\PP^1$
which for general $\tau_1$, $\tau_3$ is birational onto its image.

\pf The only nontrivial assertion is the last one. We shall show,
laboriously, that $\phi_{|E_1|}$ is birational on $C$ and hence on
every general translate of $C$. Given this, $\psi$ must either be
itself birational or be $2$-to-$1$, identifying the two components
of a general fibre of $\phi_{|D_1|}$. In that case we consider the
corresponding birational involution $\iota:A\to A$, which is biregular
because $A$ is minimal. It preserves the fibres of
$\phi_{|D_1|}$ and in particular it preserves the four double fibres
$C_1,\ldots,C_4$ which correspond to the branch points of
$C'\to\PP^1$. We know that $\iota$ is not the Kummer involution of $A$
(with some choice of origin) because in that case the $\pm
1$-eigenspaces of $\iota$ in $H^0(\cO_A(E_1))$ both have positive
dimension, so $\phi_{|E_1|}$ does not factor through~$\iota$. We also
know that $\iota$ is not translation by a $2$-torsion point of~$A$,
because in order to preserve the fibres it would have to be a
$2$-torsion point of~$C$ and then $\iota$ would preserve every translate
of~$C$ instead of interchanging different components of the general
fibre. So $\iota$ has fixed points, and they all lie on the double fibres. 

If these fixed points are isolated then $\iota$ is after all conjugate
to the Kummer involution. The alternative is that $\iota$ fixes each
$C_i$ pointwise. If we assume, as we may do, that $\cO_A(E_1)$ is a
symmetric line bundle then $\phi_{|E_1|}$ becomes equivariant for the
action of the extended Heisenberg group $H(3)^e$, as described in
[5], and in particular the ramification curve $R\contin\PP^2$ is
$H(3)^e$-invariant. But $R$ certainly includes the image of the branch
locus of~$\psi$ which in this case is the image of $\sum C_i$. Each of
these curves is of degree~$4$, so $R=R'+\sum\phi_{|E_1|}(C_i)$, and
$\deg R = 18$ so $R'$ must be of degree~$2$. It is the only
reduced degree~$2$ component of $R$, so it must be
$H(3)^e$-invariant. But it is easy to see, using the generators of
$H(3)^e$ given in~[5] that no such conic exists. (See the remark below
for an alternative argument.)

It remains to show that for general $\tau_1$, $\tau_3$, the map
$\phi_{|E_1|}|_C:C\to\PP^2$ is birational.

The linear system $|\cO_C(E_1)|$ embeds $C$ as the intersection of two
quadrics in $\PP^3=\PP H^0\big(\cO_C(E_1)\big)^*$. The restriction map
$\cO_A(E_1)\to\cO_C(E_1)$ induces
$$
0\To H^0\big(\cO_A(E_1-C)\big)\To H^0\big(\cO_A(E_1)\big)\To
H^0\big(\cO_C(E_1)\big)
$$
and $(E_1-C)^2=-2<0$ so the right-hand map is injective. So the image
of $C$ under $\phi_{|E_1|}$ is the projection of $C\contin\PP^3$ to
some $\PP^2$ which is determined by the $3$-dimensional
subspace~$H^0\big(\cO_A(E_1)\big)$.

Projection from a point $P\in\PP^3$ will map $C$ onto a double conic
if and only if $P$ is the vertex of a quadric cone
containing~$C$. Since $C$ is the intersection of two quadrics,
$h^0\big({\cal I}_{C/\PP^3}(2)\big)=2$, so there are only finitely
many (actually four) quadric cones containing~$C$ and therefore only
finitely many projections that fail to be birational on~$C$. Fixing
$\tau_1$ and letting $\tau_3$ vary we get a family of projections: if
we can show that this family is nonconstant (for some choice of
$\tau_1$) we shall have finished.

$H^0\big(\cO_A(E_1)\big)$ is spanned by the classical theta functions
$$
\theta{{\scriptstyle{0}\;\scriptstyle{{j}\over{3}}}\atopwithdelims[]
{\scriptstyle{0}\;\scriptstyle{0}}}
\left(z_1,z_2,
{
{\scriptstyle{4\tau_1}\;\scriptstyle{3\tau_1}}
\choose
{\scriptstyle{3\tau_1}\;\scriptstyle{\phantom{3}\tau_3}}
}
\right)
= 
\Sum_{m,n\in\ZZ}e^{\left\{\pi\sqrt{-1} (m,n+{{j}\over{3}}) 
{
{\scriptstyle{4\tau_1}\;\scriptstyle{3\tau_1}}
\choose
{\scriptstyle{3\tau_1}\;\scriptstyle{\phantom{3}\tau_3}}
}
{\scriptstyle{m}\choose{\scriptstyle{n+{{j}\over{3}}}}}
+2\pi\sqrt{-1}(mz_1+(n+{{j}\over{3}})z_2)\right\}}
$$
where we have chosen $E_1$ so that $\cO_A(E_1)$ has characteristic zero
with respect to the decomposition determined by the period
matrix~$\Pi$. We use [15] as our general reference for this theory. 

The restriction of this bundle to~$C$ is of
characteristic zero with respect to the decomposition
$\ZZ\oplus\ZZ\tau_1$ since $\gamma(1)$ and $\gamma(\tau_1)$ are in the
sublattices $\ZZ\bdf_3+\ZZ\bdf_4$ and $\ZZ\bdf_1+\ZZ\bdf_2$
respectively. In particular if we fix $\tau_1$ the bundle $\cO_C(E_1)$
does not depend on $\tau_3$. If we restrict these theta functions
to~$C$ we shall get (non-classical) theta functions determining a
$3$-dimensional subspace of $H^0\big(\cO_C(E_1)\big)$ which we must
show really does vary with $\tau_3$.

We denote by $\vartheta_j(z,\tau_1,\tau_3)$ the restriction of 
$\theta{{\scriptstyle{0}\;\scriptstyle{{j}\over{3}}}\atopwithdelims[]
{\scriptstyle{0}\;\scriptstyle{0}}}$
to $\tilde C=\{(4z,3z)\mid z\in\CC\}$.
$$
\eqalign{
\vartheta_j(z,\tau_1,\tau_3)&=\Sum_{m,n\in\ZZ} e^{
\pi\sqrt{-1}[(4m^2+6mn+2mj)\tau_1+(n+j/3)^2\tau_3]}
e^{2\pi\sqrt{-1}(4m+3n+j)z}\cr
&=\sum_{n\in\ZZ}s^{(3n+j)^2}\sum_{m\in\ZZ}t^{4m^2+6mn+2mj}
e^{2\pi\sqrt{-1}(4m+3n+j)z}\cr
}
$$
where we have set $s=e^{\pi\sqrt{-1}\tau_3/3}$ and
$t=e^{\pi\sqrt{-1}\tau_1}$.

We now need some coordinates in $\PP H^0\big(\cO_C(E_1)\big)^*$. This
we do by selecting four arbitrary fixed points $z_0$, $z_1$, $z_2$,
$z_3$ on~$C$ and taking the evaluation maps at those points as a
basis. We cannot take $z_i$ to be the $2$-torsion points, however, as
that does not give a basis, since the $2$-torsion points are coplanar
in~$\PP^3$ in this embedding. Instead we pick $z_0=0$, $z_1=1/2$,
$z_2=\tau_1/2$ and $z_3=1/3$ in~$\tilde C$. Then we consider the
matrix $\Theta=\bigg(\vartheta_j(z_i)\bigg)$, $0\le j\le 2$, $0\le
i\le 3$, and its four $3\times 3$~minors
$\hat\Theta_k=\det\bigg(\big(\vartheta_j(z_i)\big)_{i\not=k}\bigg)$.
The point $(\hat\Theta_0:\hat\Theta_1:\hat\Theta_2:\hat\Theta_3)\in
\PP^3$ is the point of $\PP\big(\bigwedge\nolimits^3
H^0\big(\cO_C(E_1)\big)\big)=\PP H^0\big(\cO_C(E_1)\big)^*$ which is
the vertex of the projection induced by $\cO_A(E_1)\to\cO_C(E_1)$.
 
Now one calculates directly, writing out the first few terms of each
$\vartheta_j(z_i)$ as a power series in $s$, whose coefficients are
Laurent series (with bounded negative degree) in~$t$. From this one
can calculate 
$$
\hat\Theta_k=s^2g_{k2}(t)+s^5g_{k5}(t)+O(s^8)
$$
and then the point will depend on $t$ unless (inter alia)
$g_{02}g_{15}-g_{12}g_{05}\equiv 0$. But this can be calculated from
the Laurent expansions of $g_{kl}(t)$. I did this using MAPLE (it is
not beyond the capacity of a determined human) and found that this
expression has the constant term~$36$. As this is not zero, we are
done.~\qed

{\it Remark.\/} In fact a general abelian surface $A$ in this family has no
order~$2$ automorphisms apart from~$-1$, because the family
corresponds to an Humbert surface of discriminant~$16$ in the moduli
space of $(1,3)$-polarised abelian surfaces. The abelian surfaces that
do have extra automorphisms of order~$2$ are the product surfaces and
the bielliptic abelian surfaces, and those are parametrised by Humbert
surfaces of discriminants~$1$ and~$4$ respectively, as is shown
in~[10]. Since it is easy to see that $\psi$ has degree at most~$2$ we
could use this fact to replace the argument above if we knew
that~$\iota$ could not be translation by a $2$-torsion point in~$C$.

\prop 3.4 Given $A$, $E_1$ and $D_1$ and homogeneous coordinates
$(x_1:x_2)$ in $|D_1|$ and $(y_1:y_2:y_3)$ in $|E_1|$, we can specify a
morphism $\phi:A\to X$ by choosing a curve $D_3\in|E_1+D_1|$ and a
trivialisation $\cO\mapsim\cO(E_1+D_1-D_3)$. There is a rational map
$\pi:X\rat\PP^2\cross\PP^1$ such that $\pi\phi=\psi$, and in
particular $\pi|_{\phi(A)}$ is a morphism.

\pf According to Theorem~2.1 we need to specify $D_3$ as in the
statement of the proposition and also a homomorphism
$$
\varepsilon:M\To H^0\big(A\setminus\big(\bigcup D_i\cup\bigcup E_j\big),
\cO_A\big)^\times\contin K(A)^\times.
$$
Obviously it is enough to specify $\varepsilon$ on a basis of
$M\imic\ZZ^4$. In $X$ we have $D_1-D_2=\div\big(\bde(1,-1,0,0)\big)$,
$E_1-E_2=\div\big(\bde(0,0,1,-1)\big)$ and
$E_1-E_3=\div\big(\bde(0,0,1,0)\big)$, so we
should define $\varepsilon$ on the space $m_1+m_2=0$ spanned by these
three by putting $\varepsilon(1,-1,0,0)=x_1/x_2$,
$\varepsilon(0,0,1,-1)=y_1/y_2$ and $\varepsilon(0,0,1,0)=y_1/y_3$. We
can think of these as functions on $A$ by composing with $\phi_{D_1}$
or $\phi_{E_1}$. The trivialisation $\cO\mapsim\cO(E_1+D_1-D_3)$ then
determines $\varepsilon(1,0,0,0)$, since
$\div\big(\bde(1,0,0,0)\big)=E_1+D_1-D_3$ in~$X$.

The rational map $\pi$ is given by the projection $\cO\oplus\cO(1)
\oplus\cO(1)\to\cO(1)\oplus\cO(1)$, which is evidently equal to
$\psi\phi^{-1}$ on $\phi(A)$. The data that determine $\phi$
include data that determine $\psi$, namely $E_1$, $E_2$, $E_3$,
$D_1$, $D_2$ and $\varepsilon|_{\{m_1+m_2=0\}}$, so $\pi|_{\phi(A)}$
is well-defined and therefore a morphism.~\qed

\corollary 3.5 $\phi:A\to X$ is birational onto its image.~\qed

Next we collect some information about the singularities of $\bar
A:=\psi(A)\contin \PP^2\cross\PP^1$, for a general choice of
$\psi$. We do not need all of this information but it also clarifies
the geometric picture.

The generic $D_\eta\in|D_1|$ is a union of two smooth curves of
genus~$1$ in $A$, say $D_\eta=D^+\coprod D^-$. The linear system $E_1$
has degree~$4$ on each of these, so the fibre $\bar
A\cap\pr_1^{-1}(\eta)$ consists, for generic $\eta\in\PP^1$, of two
plane quartic curves $\bar D^+$, $\bar D^-$ with $p_g=1$. The curve
$\bar D_\eta=\bar D^+\cup\bar D^-$ has only ordinary double point
singularities: there are 20 of these, of which 16 are the points of
$\bar D^+\cap\bar D^-$ and 4 are singularities of $\bar D^\pm$ (two
on each curve). At all of these points, $\bar A$ also has a
(non-isolated) singularity. Taking the closure we get a curve
$\bar\Gamma_{\rm int}\cup\bar\Gamma_{\rm node}\contin\Sing \bar
A$, where $\bar\Gamma_{\rm int}$ corresponds to the 16 intersection
points and $\bar\Gamma_{\rm node}$ to the 4 other nodes. Take
$\bar\Gamma$ to be the union of all dimension~$1$ components of $\Sing
\bar A$. In fact $\bar\Gamma=\bar\Gamma_{\rm int}\cup\bar\Gamma_{\rm
node}$ but we shall not need this fact.

As a scheme $\Sing \bar A$ consists of $\bar\Gamma$ and perhaps some
points (possibly infinitely near to one another, possibly infinitely
near to points of $\bar\Gamma$). We shall see shortly that such points
may in practice be ignored. Denote by $\Gamma_i$ the reduced curve in
$A$ whose image in $\bar A$ is an irreducible component $\bar\Gamma_i$
of $\bar\Gamma$. The map $\psi:A\to\bar A$ fails to be an embedding
along $\Gamma_i$; in fact it maps $\Gamma_i$ 2-to-1 onto
$\bar\Gamma_i$. We need to check that the additional information
carried by $\phi$ is sufficient to separate a general pair of points
of this kind, in other words, that $\phi|_{\Gamma_i}$ is
birational. Then we shall have to deal with the $0$-dimensional part
of the singular locus that remains, but it will turn out that this is
empty.

For all of this the essential observation is the following.

\prop 3.6 For a general $A$ with period matrix as above, the line
bundle $\cO_A(E_1+D_1)$ is very ample.

\pf We have $[E_1+D_1]=\bda$ and calculating intersection numbers on
$A$ gives $(E_1+D_1)^2=\lambda'=22$; so $E_1+D_1$ determines a
polarisation of type $(1,11)$. According to Reider's theorem, in the
form of [15], 10.4.1, such a polarisation is very ample unless either
$(A,\bda)$ is a product of elliptic curves with a product polarisation
or $A$ contains an elliptic curve $J$ such that $J.(E_1+D_1)=2$.

Suppose first that a general $A$ is a product. Then there are elliptic
curves $J,J'\contin A$ such that $A\imic J\cross J'$ and
$\bda=[\cL\boxtimes\cL']$, where $\cL$ and $\cL'$ are line bundles on
$J$ and $J'$ of degrees 1 and 11 respectively. We have the
intersection numbers $J.\bda=1$, $J'.\bda=11$, $J.J'=1$. Since we are
considering a general $A$ in the surface in the moduli space given by
the condition of 3.2, we may assume that $\rho(A)=2$, so that
$\NS(A)\tens\QQ=\QQ\bda\oplus\QQ\bdb$. Suppose
$[J]=\xi\bda+\zeta\bdb$: then
$$
0=(\xi\bda+\zeta\bdb)^2_A=2(\xi+\zeta)(11\xi+3\zeta),
$$
so $[J]=\xi(\bda-\bdb)$ or $[J]=\xi(11\bda-3\bdb)$; similarly
$[J']=\xi'(\bda-\bdb)$ or $[J']=\xi'(11\bda-3\bdb)$, with
$\xi,\xi'\in\QQ$. If $[J]=\xi(\bda-\bdb)$ then $1=J.\bda
=\xi(\bda-\bdb).\bda=8\xi$, so $\xi=1/8$;
$1=J.J'=\xi\xi'(\bda-\bdb)(3\bda-11\bdb)=49\xi'/2$ (we cannot have
$[J']=\xi'(\bda-\bdb)$ in this case as then $J.J'=0$); and finally
$11=J'.\bda=\xi'(3\bda-11\bdb).\bda
={{2}\over{49}}(3\bda^2-11\bda\bdb) =-{{176}\over{49}}$, which is
absurd. If $[J]=\xi(11\bda-3\bdb)$ a similar calculation leads to the
same result.

Suppose then that a general $A$ contains an elliptic curve $J$ with
$J.(E_1+D_1)=2$. Since $\cO_A(E_1)$ is ample this implies either
$J.E_1=J.D_1=1$ or $J.E_1=2$, $J.D_1=0$. Again we may suppose
$\rho(A)=2$ and $[J]=\xi\bda+\zeta\bdb$; as above, this implies
$[J]=\xi(\bda-\bdb)$ or $[J]=\xi(11\bda-3\bdb)$. If $J.D_1=1$ then
$[J]=\xi(11\bda-3\bdb)$ and $1=2\xi(11\bda-3\bdb)(\bda-\bdb)=392\xi$
so $\xi={{1}\over{392}}$. But then
$2=J.E_1={{1}\over{392}}(3\bda-11\bdb)\bdb =-{{9}\over{196}}$ so this
is impossible.

It remains to exclude the possibility that $J.E_1=2$. If this happens
then $2[J]=[C]\in H^2(X;\ZZ)$. But we saw earlier that
$$
[C]=\bdf_1^*\wedge 4\bdf_3^*+\bdf_1^*\wedge\bdf_4^*
\in\bigwedge\nolimits^2\Hom(\Lambda,\ZZ)\imic H^2(A,\ZZ)
$$
and this is not divisible.~\qed

{\it Remark.\/} Again one could argue, less directly, that the surfaces
for which a $(1,11)$-polarisation is not very ample are the product
surfaces and the bielliptic abelian surfaces, and that those are
parametrised by Humbert surfaces different from the one that occurs here.

\corollary 3.7 For a general choice of $D_3\in|E_1+D_1|$ and a
trivialisation of $\cO_A(E_1+D_1-D_3)$, the image of the associated
map $\phi:A\to X$ has only isolated singularities.

\pf Choose an irreducible component $\bar\Gamma_i$ of $\bar\Gamma$ and
a point $P\in\bar\Gamma_i$. For general $P$, there are precisely two
points $P_1$, $P_2\in A$ such that $\psi(P_1)=\psi(P_2)$. By 3.6, the
subspace of $H^0(\cO_A(E_1+D_1))$ given by the condition
$\sigma(P_1)=\sigma(P_2)$ is proper, so there is a non-empty
Zariski-open subset $U_P\contin H^0(\cO_A(E_1+D_1))$ for which
$\sigma(P_1)\not=\sigma(P_2)$. Furthermore, given $\sigma\in
H^0(\cO_A(E_1+D_1))$, the set of points of $\bar\Gamma_i$ whose two
preimages under $\psi$ are separated by $\sigma$ is Zariski-open. If
$\sigma\in U_P$ then this open set is non-empty, and doing this for
each component and taking $\sigma$ to be in the intersection of the
$U_P$s we can find a $\sigma$ which separates the preimages of all but
finitely many points of $\bar\Gamma$. (In principle $\psi$ might kill
a tangent direction at a general point of $\Gamma_i$ instead of
identifying two distinct points. If so, the points $P_1$ and $P_2$
will be infinitely near but this makes no difference.)

Now take $D_3$ to be the set $\{\sigma=0\}$, which we may assume to be
reduced and irreducible if we like, and take $\tau\in
H^0(\cO_A(E_1+D_1))$ such that $\{\tau=0\}=D_1\cup E_1$, so that
$\tau\in H^0(\cO_A(D_1))\tens H^0(\cO_A(E_1))$. Then we take the
trivialisation of $\cO_A(E_1+D_1-D_3)$ given by $\tau/\sigma$. Now we
have enough data to separate $P_1$ and $P_2$, in other words
$\phi|_\Gamma:\Gamma\to X$ is birational onto its image. As $\phi$
is birational outside $\Gamma$, except perhaps at finitely many
points, we are done.~\qed

So the failure of $\phi$ to be birational can only be caused by its
identifying finitely many pairs of points of $A$, or killing a tangent
direction at finitely many points. The set of such pairs, respectively
tangent directions, is the set of closed, respectively embedded,
points of the double-point scheme $\tilde D(\phi)$. So it is enough to
show that $\tilde D(\phi)$ is empty; then $\phi$ will be an embedding.

\prop 3.8 For general $D_3\in |E_1+D_1|$ and trivialisation of
$\cO_A(E_1+D_1)$, the double point scheme $\tilde D(\phi)$ is empty.

\pf If $\tilde D(\phi)\not=\emptyset$, then $\codim \tilde D(\phi)=2$ by
3.7, so (see [7], p.166 for the notation and general facts)
$\tilde\DD(\phi)=[\tilde D(\phi)]$ and is a nonzero element of
$A_0(\tilde D(\phi))$, so $\DD(\phi)\in A_0(A)\imic\ZZ$ is also
nonzero. So we want to show that in fact $\DD(\phi)=0$. By the double
point formula ([7], Theorem 9.3)
$$
\eqalign{
\DD(\phi)&=\phi^*\phi_*[A]-\big(c(\phi^*T_X)c(T_A)^{-1}\big)_2\cap[A]\cr
 &=\phi^*\phi_*[A]-c_2(\phi^*T_X)\cap[A]\cr
 &=\big([\phi(A)]-c_2(T_X)\big)\cap[A]\cr
 &=0\cr}
$$
since $\big([\phi(A)]-c_2(T_X)\big).[A]=0$ in $A_0(X)$, by the choice
of the class of $\phi(A)$.~\qed

This concludes the proof of Theorem~3.1.

\startsection 4. Further remarks

We can use the abelian surfaces constructed in the previous section to
give some rank~$2$ vector bundles on
$X=\PP\big(\cO_{\PP^2}\oplus\cO_{\PP^2}\oplus
\cO_{\PP^2}(1)\oplus\cO_{\PP^2}(1)\big)$ via Serre's
construction, extending the normal bundle $\cN_{A/X}$ to the whole
of~$X$. One needs to check that $\det \cN$ is the restriction of a
line bundle $\cE$ on~$X$ with $H^1(\cE)=H^2(\cE)=0$, but this is
immediate as $\cE=K_X$. In fact each $A\contin X$ produces a rank~$2$
bundle in this way and there are many questions that might be asked
about them. For instance, are they all isomorphic? Are they
indecomposible (presumably yes)?  Can one calculate their cohomology?
Some of these questions are answered in [14] in the case of
$X=\PP^1\cross\PP^3$, where an extension of the normal bundle exists
for the same reasons.

Another series of questions raised by this example is the possibility
of extending the procedure to other~$X$. The results of Section~1
allow one to generate other possibilities among the smooth toric
$4$-folds with $\rho=2$, but the proofs in Section~3 used some special
geometry and in particular the fact that the linear system spanned by
$E_1$, $E_2$ and $E_3$ is complete. In other cases one would
presumably have to work with very far from complete linear systems and
the methods of this paper might not be adequate. In any case a
more interesting problem might be to revert to $\rho(X)=1$ but allow
singular toric varieties, and try to apply the results of Kajiwara
from~[11]. The case of weighted projective spaces is a natural
starting point. Another possibility would be to work with Batyrev's list
[2] of toric Fano $4$-folds.

The method used to prove that the morphism $A\to X$ we produce is an
embedding is very clumsy. In the case of $X=\PP^n$ one has an elegant
criterion in the form of Reider's theorem. It would be interesting to
have a way of distinguishing the embeddings (or even the birational
morphisms) among the morphisms into toric varieties, say in terms of
Cox's description in [6].

\startsection{References}

\noindent[1] V.V.~Batyrev, {\it On the classification of smooth
projective toric varieties}, T\^ohoku Math. J. {\bf 43} (1991),
569--585.

\noindent[2] V.V.~Batyrev, {\it On the classification of toric Fano
4-folds}, Preprint 1998 (math/9801107).

\noindent[3] T.~Bauer \& T.~Szemberg {\it Abelian threefolds in
$(\PP_2)^3$}, in {\it Abelian Varieties} (W.~Barth, K.~Hulek \&
H.~Lange, Eds.) 19--23, De Gruyter, Berlin 1995.

\noindent[4] Ch.~Birkenhake, {\it Abelian threefolds in products of
projective spaces}, Abh. Math. Sem. Univ. Hamburg {\bf 65} (1995),
113--121.

\noindent[5] Ch.~Birkenhake \& H.~Lange, {\it A family of Abelian
surfaces and curves of genus four}, Manuscr. Math. {\bf 85} (1994),
393--407.

\noindent[6] D.A.~Cox, {\it The functor of a smooth toric variety},
T\^ohoku Math.\ J., {\bf 47} (1995), 251--262.

\noindent[7] W.~Fulton, {\it Intersection theory}, Springer, New York (1984).

\noindent[8] M.A.~Guest, {\it The topology of the space of rational
curves on a toric variety}, Acta Math. {\bf 174} (1995), 119-145.

\noindent[9] K.~Hulek, {\it Abelian surfaces in products of projective
spaces}, in {\it Algebraic Geometry, L'Aquila} (A.J.~Sommese,
A.~Biancofiore \& E.L.~Livorni, eds.) LNM~{\bf 1417}, 129--137,
Springer, Berlin 1988.

\noindent[10] K.~Hulek \& S.~Weintraub, {\it Bielliptic abelian surfaces},
Math. Ann. {\bf 283} (1989), 411--429.

\noindent[11] T.~Kajiwara, {\it The functor of a toric variety with
enough invariant effective Cartier divisors}, T\^ohoku Math. J. {\bf
50} (1998), 139--157.

\noindent[12] P.~Kleinschmidt, {\it A classification of toric varieties
with few generators}, Aequationes Math. {\bf 35} (1988), 254--266.

\noindent[13] H.~Lange, {\it Abelian surfaces in $\PP^1\cross\PP^3$},
Arch. Math. {\bf 63} (1994), 80--84.

\noindent[14] H.~Lange, {\it A vector bundle of rank~$2$ on
$\PP^1\cross\PP^3$}, J. London Math. Soc. {\bf 57} (1998), 583--598.

\noindent[15] H.~Lange \& Ch.~Birkenhake {\it Complex abelian
varieties}, Springer, Berlin 1992.

\noindent[16] T.~Oda, {\it Convex bodies and algebraic geometry},
Springer, Berlin (1987).

\noindent[17] A.~Van~de~Ven, {\it On the embedding of abelian varieties
in projective space}, Annali di Matematica Pura ed Applicata {\bf 103}
(1971), 127--129.
 
\end